\documentclass[12pt]{article} 
\usepackage{amsmath, amssymb, theorem, latexsym}
\def\nz{\ifmmode {I\hskip -3pt N} \else {\hbox {$I\hskip -3pt N$}}\fi}
\def\zz{\ifmmode {Z\hskip -4.8pt Z} \else
       {\hbox {$Z\hskip -4.8pt Z$}}\fi}
\def\qz{\ifmmode {Q\hskip -5.0pt\vrule height6.0pt depth 0pt
       \hskip 6pt} \'else {\hbox
       {$Q\hskip -5.0pt\vrule height6.0pt depth 0pt\hskip 6pt$}}\fi}
\def\rz{\ifmmode {I\hskip -3pt R} \else {\hbox {$I\hskip -3pt R$}}\fi}
\def\cz{\ifmmode {C\hskip -4.8pt\vrule height5.8pt\hskip 6.3pt} \else
       {\hbox {$C\hskip -4.8pt\vrule height5.8pt\hskip 6.3pt$}}\fi}

\def\tr {{\rm \;Tr \;}}
\def\Tr {{\rm \;Tr \;}}
\def\Op {{\rm \;Op \;}}
\def\Opw {{\rm \;Op^w \;}}

\newcommand {\pa}{\partial}

\makeatletter
\@addtoreset{equation}{section}

\makeatother

\newtheorem{theorem}{Theorem}[section]
\newtheorem{lemma}[theorem]{Lemma}
\newtheorem{proposition}[theorem]{Proposition}

\newtheorem{remark}[theorem]{Remark}
\newtheorem{corollary}[theorem]{Corollary}

\title{Non linear eigenvalues and analytic hypoellipticity.}
\author{Sagun  Chanillo, \\Department of Mathematics\\
Rutgers, The State University of New Jersey,\\
Hill Center for the Mathematical Sciences, Busch Campus,\\
 Piscataway, New Jersey 08854, USA\\
Bernard Helffer,\\
D\'epartement de Math\'ematiques, UMR CNRS 8628,\\
Universit\'e Paris-Sud, Bat. 425 \\
91405 Orsay Cedex,  France\\
and\\
 Ari  Laptev,\\
Department of Mathematics,\\
Kungl. Tekniska H\"ogskolan,\\
10044 Stockholm, Sweden.
} 
\begin{document}
\bibliographystyle{plain}
\maketitle
\begin{abstract}
Motivated by the problem of analytic hypoellipticity, we show 
 that  a special family  of compact non selfadjoint operators
 has a non zero eigenvalue. We recover old results obtained
 by ordinary differential equations techniques and show how
it can be applied to the higher dimensional case. This gives in particular a new class of hypoelliptic, but  not analytic hypoelliptic operators.
\end{abstract}
\section{Introduction}

There is a long history highlighting the links between spectral analysis and the construction of hypoelliptic but not analytic hypoelliptic operators. Since the basic works  of \cite{Met1, Tr1, Tar, Sj, GS} and the necessary conditions obtained by \cite{Met2}, there has been a lot of effort in 
 understanding when H\"ormander sums of squares operators
 formed by real-analytic vector fields fail to satisfy the analytic hypoellipticity property. These results more or less may be summarized by the fact that failure of analytic hypoellipticity
 occurs whenever the characteristic set of the vector fields satisfies a certain condition conjectured by Tr\`eves \cite{Tr2}.\\
 We refer to \cite{BaGo, Cha,Chr1,Chr2,Chr3,Chr4,Chr5,Co, HH1,HH2,HH3,
He0, He1, He2,Ho,Ol, OlRa, PhRo} for various examples. Two types of problems
 appear. The first type is described by the Baouendi-Goulaouic example \cite{BaGo}.
For showing that $D_{x_1}^2 + x_1^2 D_{x_2}^2 + D_{x_3}^2$ is not
 hypoanalytic, it
 is shown that it is enough to find a complex $\lambda$
 such that $D_{x_1}^2 + x_1^2 + \lambda^2$ is not injective. It is enough to take $\lambda = i \sqrt{\lambda_j}$
 where $\lambda_j$ is an eigenvalue of the harmonic oscillator. This idea can be used in a quite general context, see
 \cite{He3} and \cite{Cha} for more recent variants, without any restrictions on the dimension.\\
The second type was initially proposed by B. Helffer in \cite{He1, He2}
 and solved by Pham The Lai-Robert \cite{PhRo}.
For showing that the operator $D_{x_1}^2 + (x_1^2 D_{x_2}- D_{x_3})^2$
 is not analytic hypoelliptic, one has to show that
  it is enough to find a complex $\lambda$
 such that $D_{x_1}^2 + (x_1^2 -\lambda)^2$ is not injective.
 This problem is more involved. The proof in  \cite{PhRo}
 although multi-dimensional in principle seems to break down  almost immediately
 when the spectral problem is in dimension greater than $1$. The conditions of Theorem 2.3 in  \cite{PhRo} (Section 3, Application 1) are not so easy to verify. On the other hand, these authors  prove the existence of a complete system
 of eigenvectors. This property is much stronger but not useful for the problem
 of non analytic hypoellipticity, which requires only the existence of one eigenvector.
After this work,  M. Christ (and then many others as recalled in  the references above)
 extended this example. Typically M. Christ can deal with the family
$D_{x_1}^2 + (x_1^m -\lambda)^2$ ($m >1$), in particular with $m$ odd which seems not accessible by the Pham The Lai-Robert method \cite{PhRo} \cite{Rob2}. 

The method of Christ relies on the Wronskian function and thus seems limited to models which give rise to one dimensional spectral problems. Our aim is to propose a technique
 permitting to treat many new examples not necessary in dimension $1$.

Our family of operators would be of the type
\begin{equation}
H(x, D_x, \lambda) = - \Delta + ( \lambda - P(x))^2\;,
\end{equation}
where $x\mapsto P(x)$ is an homogeneous elliptic polynomial on $\mathbb  R^n$
 of order $m>1$. \\
Although it could be a rather natural conjecture that in this case
 there exists always  $\lambda \in \mathbb C$ such that
 $H(x, D_x, \lambda)$ is non injective on $\mathcal  S(\mathbb R^n)$,
 our results will be only true  for $n \leq 3$ and $ m \geq m(n)>1$ (See Theorems \ref{Propositionn=2} and 
 \ref{newspect}).

The spectral result which is considered can first be reduced
 to a problem for a compact operator.

We rewrite $H(x, D_x, \lambda)$ in the form
\begin{equation}
H(x, D_x, \lambda) = L - 2 \lambda M + \lambda^2\;,
\end{equation}
with 
\begin{equation}
L = - \Delta + P(x)^2\;,\; M = P(x)\;.
\end{equation}
 The operator $L$ 
 is invertible and its inverse is a pseudo-differential operator (See
  appendix \ref{appendixC} and Helffer \cite{He4}).
It is also easy to give sufficient condition for determining whether 
 the operator 
\begin{equation}
A:= L^{-1}
\end{equation}
 belongs to a given Schatten class (see \cite{Ron} and  appendix \ref{appendixB}). The Hilbert-Schmidt
 character can be deduced from the fact that the Weyl symbol is in $L^2( \mathbb R^n\times \mathbb R^n)$.
 The restriction $n \leq 3$ appears for example if $m \geq 2$
 and if we want to have $A:= L^{-1}$ Hilbert-Schmidt. The condition
 that $A$ is Trace class leads to $m > 1$ and $n=1$.

Then the initial problem is reduced to the spectral analysis of 
\begin{equation}\label{pbnl0}
(I - 2 \lambda B + \lambda^2 A) u = 0\;.
\end{equation}
with 
\begin{equation}
B = A^\frac 12 P A^\frac 12
\end{equation}
In the spirit of \cite{PhRo}, one is led to the study of 
 the so-called operator pencils for which there is a large literature, for e.g. Markus's  book  \cite{Ma}. Additional literature was mentioned to us by Markus.  However these results do not apply to our situation. Typically  one has results where the operator pencils are of the type
$$
I - 2 \lambda B - \lambda^2 A\;,
$$
where $A,B$ are selfadjoint and compact, see Friedman-Shinbrot \cite{FrSh} and reference therein. Our situation is what is called
 in the literature an elliptic pencil.

A few months ago, one of the authors (S.C.) proved a result  \cite{Cha1}, which we later realized  was 
 a weak version  of   Lidskii's Theorem. Motivated by \cite{Cha1}, 
 we were led to consider the computation of traces in the spectral problems we will deal with in this article.  Lidskii's Theorem will  systematically be applied in the sequel.

\paragraph{Acknowledgements}.\\
We thank A.S.~Markus, D.~Robert and M.~Solomyak for useful correspondence. S.C.  wishes
 to thank F. Tr\`eves for his encouragement and Shri S. Devananda for useful comments. B. H. and A.L. thank the Mittag-Leffler Institute and a partial support by the SPECT ESF european programme. The research of S.C. was supported
 in part by a grant from the NSF.

\section{Lidskii's Theorem and applications}
Let us  show how to use  Lidskii's Theorem.
We consider the problem of determining if there exists a non trivial
 pair $(\lambda,v)$ such that
\begin{equation}\label{pbnl}
(I - 2 \lambda B + \lambda^2 A) u = 0\;.
\end{equation}

The initial motivating example is the  example
 where~:
\begin{equation}
L = D_t^2 + t^{2m}\;,\; A = L^{-1}\;,\; B = A^\frac 12 t^m A^\frac 12
\end{equation}
which was solved by Pham The Lai-Robert \cite{PhRo}, when $m >0 $ is even
 and by Christ \cite{Chr1} when $m >1$ is odd.

We first use the reduction to the linear spectral problem.
It is enough to show that the operator $\mathcal D$ defined
 by
\begin{equation}
\mathcal D := \left( \begin{array}{cc}
2 B & A^\frac 12 \\
- A^\frac 12& 0 
\end{array}
\right )
\end{equation}
has a non zero eigenvalue $\mu$. The first component
 of the eigenvector is an eigenvector of the problem (\ref{pbnl})
 with $\mu = \frac 1 \lambda$.\\

If $B$ and $A$ are compact, $\mathcal D$ is compact but the main difficulty is
 that $\mathcal D$ is not selfadjoint. Standard results as for example explained in \cite{Sim} do not apply.

We would like to use Lidskii's Theorem (see \cite{Sim} or \cite{BiSo}) in the form
\begin{theorem}.\\
Let $\mathcal C$ be a trace class operator then 
$$\sum_j \lambda_j (\mathcal C) = \Tr \mathcal C\;.
$$
In particular, if the spectrum $\sigma(\mathcal C)$ satisfies 
$$
\sigma ( \mathcal C) = \{ 0\}\;,$$
then
$$
\Tr \mathcal C^k = 0 \;,\; \forall k \in \mathbb N^*\;.
$$
\end{theorem}
As an immediate corollary, we get~:
\begin{corollary}{Rank 2 criterion}\label{Cor22}.\\
If $\mathcal D$ is Hilbert-Schmidt (that is $B$ Hilbert-Schmidt
 and $A$ positive and Trace class) and if the condition~:
$$
\Tr ( 2 B^2 - A) \neq 0\;,
$$
is satisfied, then $\mathcal D$ has at least one non zero eigenvalue.
\end{corollary}
{\bf Proof}.\\
The proof is by contradiction. If $\mathcal D$ has no non zero 
eigenvalue, the same is true for $\mathcal C = \mathcal D^2$.
We then apply the theorem to $\mathcal C$ with $k=1$.\\

One could also try to use the criterion for other values of $k$.
If we first consider the case $k=1$, one gets
 that if $A^\frac 12$ and $B$ are Trace class and if
 $\Tr B \neq 0$ then $\mathcal D$ has at least one non zero eigenvalue.
In our applications (where $A = (- \Delta + P(x)^2)^{-1}$), this  is not very useful, because  the condition on $A^\frac 12$ is too strong and never satisfied. The consideration of the cases  $k=3$ and $k=4$ will leads to interesting and new results.  One will exploit the  two following corollaries.
\begin{corollary}{Rank 3 criterion}\label{Cor23}.\\
If $A^\frac 32$ and $B^3$ are trace class, then, if 
\begin{equation}\label{kegal3}
\Tr \left( 4 B^3 - 3 B A \right) \neq 0\;.
\end{equation}
is satisfied, then $\mathcal D$ has at least one non zero eigenvalue.
\end{corollary}

\begin{corollary}{Rank 4 criterion}\label{Cor24}.\\
If $A$ and $B^2$ are  Hilbert-Schmidt, then, if 
\begin{equation}\label{kegal4}
\Tr \left( 8  B^4 - 8 B^2 A +  A^2 \right) \neq 0\;.
\end{equation}
is satisfied, then $\mathcal D$ has at least one non zero eigenvalue.
\end{corollary}

\section{Application  of the rank 2 criterion}
\subsection{The Christ--Hanges-Himonas--Pham The Lai-Robert example}
\begin{theorem}\label{ThCh}.\\
If $m>1$, the problem
$$ \left(D_t^2 + (t^m  - \lambda )^2 \right)f =0\;,
$$
has a solution $(\lambda,f)$ with $\lambda \in \mathbb C$ and $f\in \mathcal S(\mathbb R^n)$,
$f\not\equiv 0$. 
\end{theorem}
{\bf Proof}.\\
Let us show that
 the condition in Corollary \ref{Cor22} is satisfied. 
Using that $(D_t^2 + \gamma t^{2m})$ is isospectral
 to $ \gamma^{\frac{1}{m+1}} (D_t^2 + t^{2m})$, one gets first the identity
$$
 \Tr  (D_t^2 + \gamma t^{2m})^{-1} = \gamma^{ - \frac{1}{m+1}}\Tr  (D_s^2 + s^{2m})^{-1}
$$
Differentiating with respect to $\gamma$ and taking $\gamma =1$, leads
 to  
\begin{equation} \label{identiteutile}
\frac{1}{m+1} \Tr \left( (D_t^2 + t^{2m})^{-1} \right) =
 \Tr \left( (D_t^2 + t^{2m})^{-1}t^{2m} (D_t^2 + t^{2m})^{-1} \right)\;.
\end{equation}
It is indeed
 enough to see that, if $C$ is Hilbert-Schmidt, then
\begin{equation}\label{causch}
\Tr C^2 = \langle C\;,\; C^* \rangle_{H.S} \leq || C ||_{H.S}\;\cdot\;|| C^* ||_{H.S}
 = \Tr C C^*
\end{equation}
Let us see how it is used in our case.
We observe that, by cyclicity of the trace (see Proposition \ref{PropCyc}
 in  appendix \ref{appendixA}),
 we have
$$
\Tr B^2 = \Tr C^2\;,
$$
with $C = t^m(D_t^2+ t^{2m})^{-1}$, if  
$B = (D_t^2+ t^{2m})^{-\frac 12}t^m(D_t^2+ t^{2m})^{-\frac 12}$.
We then get that
$$
\begin{array}{l}
\Tr C C^* = \Tr t^m(D_t^2+ t^{2m})^{-2} t^m = \\
\Tr t^{2m}(D_t^2+ t^{2m})^{-2} = \\
\Tr (D_t^2+ t^{2m})^{-1}t^{2m}(D_t^2+ t^{2m})^{-1}
\end{array}
$$
which is the quantity which was computed in (\ref{identiteutile}). We note
 that this 
time, we do not have anymore the restriction that $m$ is even
 for applying the results.

This gives:
\begin{equation}\label{ineqfin}
\Tr \left( 2 B^2 - A\right) = (\frac{2}{m+1} - 1) \Tr (A) < 0\;,
\end{equation}
if $ m >1$.

\subsection{The Hoshiro-Costin-Costin example}
Let us now try to recover results by Hoshiro \cite{Ho} and O. and R. Costin
 \cite{Co}. The goal will be partially achieved by the
\begin{theorem}\label{ThHo}.\\
If
\begin{equation}\label{ineg}
2 \ell +1 < m\;,
\end{equation}
then 
the problem
$$ \left(D_t^2 + (t^m  - t^\ell \lambda )^2\right) f =0\;,
$$
has a solution $(\lambda,f)$ with $\lambda \in \mathbb C$ and $f\in \mathcal S(\mathbb R^n)$,
$f\not\equiv 0$. 
\end{theorem}
We expand the operator in the usual way:
\begin{equation}
\begin{array}{l}
I - 2 \lambda (D_t^2+ t^{2m})^{-\frac 12}t^{\ell+m}(D_t^2+ t^{2m})^{-\frac 12}  + \lambda^2 (D_t^2+ t^{2m})^{-\frac 12}t^{2\ell}(D_t^2+ t^{2m})^{-\frac 12}\\
 = I -2 \lambda B + \lambda^2 A\;.
\end{array}
\end{equation}
Here $$B = (D_t^2+ t^{2m})^{-\frac 12}t^{\ell+m}(D_t^2+ t^{2m})^{-\frac 12}$$
 and $$A = (D_t^2+ t^{2m})^{-\frac 12}t^{2\ell}(D_t^2+ t^{2m})^{-\frac 12}\;.$$
We note that $\ell$ should satisfy
$$
0\leq \ell < m\;.
$$
We observe that
$$
\begin{array}{l}
\tr B^2 = \tr  (D_t^2+ t^{2m})^{-\frac 12}t^{\ell+m}(D_t^2+ t^{2m})^{-1}t^{\ell+m}(D_t^2+ t^{2m})^{-\frac 12}\\
 = \tr (t^{m}(D_t^2+ t^{2m})^{-1}t^{\ell} t^m (D_t^2+ t^{2m})^{-1} t^\ell)
\end{array}
$$
We take 
$$
C= t^{m}(D_t^2+ t^{2m})^{-1}t^{\ell}\;.
$$
We get as before the estimate
\begin{equation}
\begin{array}{l}
\tr B^2 \leq \tr t^{2m}(D_t^2+ t^{2m})^{-1}t^{2\ell}(D_t^2+ t^{2m})^{-1}\\
= \tr t^{2\ell}(D_t^2+ t^{2m})^{-1}t^{2m}(D_t^2+ t^{2m})^{-1} \;.
\end{array}
\end{equation}
For computing the right hand side, we introduce as before 
a parameter $\gamma$ and observe that
$$
\Tr (  t^{2\ell} (D_t^2+ \gamma t^{2m})^{-1} ) =  \gamma^{ - \frac{\ell+1}{m+1}}
\Tr (  s^{2\ell} (D_s^2+ s^{2m})^{-1} )\;.
$$
Differentiating with respect to $\gamma$,
 we get
\begin{equation}
 \Tr (  t^{2\ell} (D_t^2+  t^{2m})^{-1} t^{2m}(D_t^2+  t^{2m})^{-1} )
=  \frac{\ell+1}{m+1} \Tr (  s^{2\ell} (D_s^2+ s^{2m})^{-1} )
\end{equation}
This finally gives
\begin{equation}
\tr (2 B^2 - A) \leq \left (2 \frac{\ell+1}{m+1} -1 \right) \tr A < 0\;,
\end{equation}
\section{The main tools}
Four tools were employed in the arguments in the preceding sections. In this section we elaborate briefly on these tools. The tools apply once the trace for our operators is defined. The necessary lemmas needed to prove the existence of the various traces which come
 up in our arguments are presented in the Appendix. The four tools we need are~:
\begin{enumerate}
\item Invariance by cyclicity of the trace,
\item Scaling invariance of $P$ and $A_\gamma$,
\item Cauchy-Schwarz inequality in the Hilbert-Schmidt spaces and positivity,
\item Invariance by taking the adjoint.
\end{enumerate}
\paragraph{Cyclicity.} The justification of the formula 
$$
\Tr (CD) = \Tr (DC)\;,
$$
where $C$ and $D$ are Hilbert-Schmidt can be extended
 slightly using the results of Appendix \ref{appendixA}. We will systematically
 identify various non commutative polynomial of $P$ and $A$ giving the same trace.
\paragraph{Scaling.}
 We introduce 
$$
A_\gamma = (- \Delta + \gamma  P^2)^{-1}\;,\; A_1 = A\;, \; B = A^\frac 12 P A^\frac 12\;.
$$
We also observe that $P$ and $A$ are selfadjoint and that $A$ is positive.
We shall also use that $P$ is homogeneous of degree $m$ with respect
 to a dilation and that $- \Delta$ is homogeneous of degree $-2$.
Under  this condition, we have immediately by dilation~:
\begin{lemma}\label{LemHom}.\\
 $A_\gamma$ is isospectral to $\gamma^{- \frac{1}{m+1}} A_1$.
\end{lemma}
As a corollary, we get, under the assumption that the objects in consideration are trace class
\begin{equation}\label{scalell}
\Tr A_\gamma^\ell = \gamma^{- \frac{\ell}{m+1}} \Tr A^\ell \;.
\end{equation}
\paragraph{Cauchy-Schwarz and positivity.}
For a pair of Hilbert-Schmidt operators $C$, $D$ we will use the properties (with some variants)~:
\begin{equation}\label{positif}
\Tr C C^* \geq 0\;,
\end{equation}
and 
\begin{equation}\label{gencausch}
\Tr C D^* \leq \sqrt{ \Tr C C^*}\sqrt{\Tr D D^*}\;.
\end{equation}
We recall that we used this with $D=C^*$ in (\ref{causch}). 
\paragraph{Invariance by taking the adjoint.}
It is  well known, that $\Tr C^* = \overline{\Tr C}$. If we observe here
 that {\bf our operators are real operators}, we also have~:
\begin{equation}
\Tr C = \Tr C^* \;.
\end{equation}

\section{Application of the rank 3 criterion}
In order to apply Corollary \ref{Cor23}, we need to verify (\ref{kegal3})
\begin{equation*}
 4 \Tr B^3 - 3 \Tr BA \neq 0
\;,
\end{equation*}
and to verify that $A^\frac 32$ and $B^3$ are trace class. We will assume in this section
 that the homogeneous polynomial $P$ is elliptic. Thus we also have without loss of generality,
\begin{equation}\label{Ppositif}
P \geq 0\;.
\end{equation}
Using the ellipticity of $P$ and  (\ref{C.3}), we easily see that $A^\frac 32$ and $B^3$ are trace class provided $n=2$, $m\geq 4$. We have
\begin{lemma}\label{Lemman=2}.\\
Assume $n=2$, $m\geq 4$ and let $P$ be a homogeneous elliptic polynomial.
Then
\begin{equation}
\Tr (4 B^3 - 3 BA) \leq \left( 2 \, \frac{m+2}{m+1} -3\right) \Tr (BA) < 0\;.
\end{equation}
\end{lemma}
{\bf Proof~:}\\
The strict inequality in the statement of Lemma \ref{Lemman=2} follows from
 the fact that $P$ is elliptic, non negative and $m\geq 4$. The conditions
 $n=2$, $m\geq 4$, ensure as noted above that the traces that occur in Lemma \ref{Lemman=2} and in the ensuing proof are all defined. Now,
\begin{equation}\label{identities}
\begin{array}{ll}
\Tr (B^3) &= \Tr (PA)^3\;,\\
\Tr (BA) &= \Tr (PA^2)\;.
\end{array}
\end{equation}
We will establish,
\begin{equation}\label{trickn=3}
\Tr (PA)^3 \leq \frac 12 \left( \frac{m+2}{m+1}\right) \Tr (PA^2)\;.
\end{equation}
Combining (\ref{trickn=3}) with (\ref{identities}) we get
\begin{equation}\label{5.4}
\Tr(B^3) \leq \frac 12 \left( \frac{m+2}{m+1}\right) \Tr (BA)\;.
\end{equation}
Our lemma follows easily from (\ref{5.4}). We now prove (\ref{trickn=3}).
The scaling argument is used in the following way~:
\begin{equation}
\Tr (P A_\gamma)^3 = \gamma^{- \frac 32 \frac{m+2}{m+1}} \Tr (P A)^3\;.
\end{equation}
By differentiation, we get
\begin{equation}\label{5.6}
\Tr \left( (PA)^3 P^2 A\right) = \frac 12 \frac{m+2}{m+1} \Tr (PA)^3\;.
\end{equation}
Since $P\geq 0$, the  Cauchy-Schwarz inequality gives~:
$$
\begin{array}{ll}
\Tr (PA)^3&= \Tr \left( (AP^\frac 12) (P^\frac 12 APAP)\right) \\
&\leq \left ( \Tr (AP^\frac 12 \; P^\frac 12 A)\right)^\frac 12 \left ( \Tr (P^\frac 12 APAP \cdot PAPA P^\frac 12)\right)^\frac 12 \\
& = \left ( \Tr (P A^2)\right)^\frac 12 \left ( \Tr (PA)^3P^2A\right)^\frac 12\;.
\end{array}
$$
Using (\ref{5.6}), we get
$$
\Tr (PA)^3 \leq \left ( \Tr (P A^2)\right)^\frac 12 \left ( \frac 12 \frac{m+2}{m+1} \Tr (PA)^3\right)^\frac 12\;.
$$
So this implies (\ref{trickn=3}).
To summarize, we have proved
\begin{theorem}\label{Propositionn=2}.\\
If $n=2$, $m\geq 4$ and if $P$ is an elliptic  positive homogeneous polynomial of degree
 $m$, then there exists a non trivial solution $(\lambda,f)$ in $\mathbb C \times \mathcal S(\mathbb R^2)$ of
$$
(- \Delta + (P(x) -\lambda)^2)f=0\;.
$$
\end{theorem}

\section{Application of the rank 4 criterion}
In this section we will use Corollary \ref{Cor24}. For the formal part of the argument
 it is not necessary to assume that $P$ is an elliptic polynomial or positive,
 in contrast to the previous section. However by assuming ellipticity on $P$, we easily
 verify using (\ref{C.3}) that $A$ is Hilbert-Schmidt and $B^4$ is trace class when,
$$
- 4 + n (1 + \frac 1m) < 0\;.
$$
This imposes a dimensional restriction, $n\leq 3$, and $m > 3$. See
 also Remark~\ref{Remark6.3}. There is no dimensional restriction in the formal part of the argument. We have,
\begin{lemma}\label{Lemma6.1}.\\
Let $n \leq 3$, $m\geq 6$ and $P$ a homogeneous elliptic polynomial of degree $m$.
 Then
\begin{equation}
\Tr \left( 8 B^4 - 8 B^2 A + A^2\right) \geq \Tr (8 B^4) + \left(\frac{m-7}{m+1}\right)
 \Tr A^2\;, \mbox{ for } m\geq 7\;,
\end{equation}
and 
\begin{equation}
\Tr \left( 8 B^4 - 8 B^2 A + A^2\right) \geq  \frac{7m-41}{8(m+1)}
 \Tr A^2\;, \mbox{ for } m\geq 6\;.
\end{equation}
\end{lemma}
{\bf Proof~:}\\
As observed above via (\ref{C.3}) the traces that occur in the statement of Lemma \ref{Lemma6.1} and the arguments to follow are all defined since $n\leq 3$ and $m\geq 5$. Our
 lemma easily follows from,
\begin{equation}\label{cascal}
\Tr (B^2 A) \leq \frac{1}{m+1} \Tr (A^2)\;,
\end{equation}
and
\begin{equation}\label{+causch}
8 \Tr (B^2 A) \leq \left( \frac{6}{m+1} + \frac 18\right) \Tr A^2 + 8 \Tr B^4\;.
\end{equation}
We begin with the proof of (\ref{cascal}). We have,
$$
\Tr B^2 A = \Tr (A^\frac 12 P A^\frac 12)^2 A = \Tr (PA)^2 A\;.
$$
We will use the Cauchy-Schwarz inequality in two different ways.
The first trivial idea is to write
\begin{equation}\label{cash4}
\Tr B^2 A \leq \frac {\alpha}{2} \Tr B^4 + \frac {1}{2 \alpha} \Tr A^2\;,
\end{equation}
which is true for any $\alpha \in ]0,1[$.\\
Using the cyclicity of the trace, this can equivalently be written in the form
\begin{equation}\label{cash4a}
\Tr (PA)^2 A\leq \frac {\alpha}{2} \Tr (PA)^4 + \frac {1}{2 \alpha} \Tr A^2\;.
\end{equation}
It is immediate to see that this inequality is not sufficient for getting
the expected  inequality
\begin{equation}
8 \Tr B^2  A <  8 \tr B^4 + \Tr A^2\;.
\end{equation}
So we try an alternative Cauchy-Schwarz inequality, by writing
$$
\begin{array}{ll}
\Tr (PA)^2 A & = \Tr A^\frac 12 P A\; P A^\frac 32 \\
&\leq ( \Tr A^\frac 12 P A A P A^\frac 12)^\frac 12 ( \Tr P A^\frac 32 A^\frac 32 P)^\frac 12\\
& \leq ( \Tr P  A^2P A)^\frac 12 ( \Tr P^2 A^3)^\frac 12\;.
\end{array}
$$
This leads to
\begin{equation}
\Tr (PA)^2 A \leq \Tr P^2 A^3\;.
\end{equation}
We now use the scaling invariance. As we have seen in (\ref{scalell}),
 we have
\begin{equation}
\Tr A_\gamma^2 = \gamma ^{ - \frac {2} {m+1} }\; \Tr A^2\;,
\end{equation}
and differentiating with respect to $\gamma$ and taking $\gamma=1$, we get
\begin{equation}
\Tr A^3 P^2 = \frac {1} {m+1} \; \Tr A^2\;,
\end{equation}
This leads to (\ref{cascal}).
We now prove (\ref{+causch}). We now 
 combine the inequalities (\ref{cash4a}) and (\ref{cascal}).
We write
$$
\begin{array}{ll}
8  \tr AB^2 & = 6 \Tr A B^2 + 2 \Tr A B^2 \\
&\leq \frac {6}{m+1} \Tr A^2 + \alpha \Tr A^2 + \frac 1 \alpha \Tr B^4\;.
\end{array}
$$
The choice of $\alpha = \frac 18$ gives (\ref{+causch}).
We leave as an exercise for  the reader that this idea cannot give
 a better condition on $m$. Collecting our results, we have shown the
\begin{theorem}\label{newspect}.\\
Let $n\leq 3$. Let $P(x)$ be a homogeneous polynomial of degree
 $m$, $m \geq 6$, which is elliptic, i.e. $P(\sigma )\neq 0$ if $\sigma \in S^{n-1}$. Then the problem
$$
- \Delta f + (P(x) - \lambda)^2 f =0\;,
$$
has a solution $(\lambda,f)$ with $f\in \mathcal S(\mathbb R^n)$,
$f\not\equiv 0$.
\end{theorem}

\begin{remark}\label{Remark6.3}.\\
The hypothesis that $P$ be elliptic can perhaps be relaxed in the  spirit of \cite{Cha}. 
For example in two dimensions, if one imposes the condition that the diameter
of the tubes $-1<P(x,y)<1$  tapers fast enough, one recaptures compactness properties (see also \cite{HeNo}). However
 one could be then forced to study higher order traces. This is because the
 $p$ value of the Schatten class  $\mathcal C_p$  to which
 the operator $L^{-1}$ belongs to will in general be large. The example
 when $n=2$ and $P(x_1,x_2)= x_1x_2 (x_1^2+ x_2^2)^k$ for $k$ large does satisfy
 the hypotheses of Corollary \ref{Cor24} and thus we obtain the conclusions of Theorem~\ref{newspect}.
\end{remark}

\section{Application to failure of analytic hypoellipticity}
Let us collect some of the standard consequences
 of our spectral analysis.
By applying Theorem \ref{ThHo}, we get
\begin{proposition}.\\
If $2 k +1 < m$,  the operator
$D_t^2 + (t^m D_y - t^k D_z)^2$ is not analytic hypoelliptic.
\end{proposition}
This recovers for $k=1$ all the mentioned known results with
 a unified  elementary  proof 
 but gives for $k>1$ only  partially results by Hoshiro \cite{Ho} and O. and R. Costin
 \cite{Co}.\\
A consequence of Theorem \ref{newspect} is the following 
\begin{proposition}.\\
The operator
$$
P_k:= \sum_{j=1}^p D_{x_j}^2  + \left( \left(\sum_{j=1}^p x_j^2\right)^{k} D_{x_{p+1}} - D_{x_{p+2}}\right)^2\;,
$$
is not  analytic hypoelliptic in the following cases~:
\begin{itemize}
\item $p=2$, $k\geq 2$\;,
\item $p=3$, $k \geq 3$\;.
\end{itemize}
\end{proposition}
{\bf Proof}.\\
The smooth solution to $P_k u=0$ that is not real-analytic can be constructed in
a
neighborhood of the origin by means of the formula,
$$ u(x,x_{p+1},x_{p+2})=\int_0^\infty\exp(i\rho^{2k+1}x_{p+1}+i\rho\lambda
x_{p+2})f(\rho x) \exp(-M\rho) d\rho,$$
where $x=(x_1,\ldots,x_p)$ and $f$ is the eigenfunction we have constructed
in Theorem \ref{newspect} and $M>0$ picked suitably large so that the integral converges
for $x_{p+2}$ in some interval centered at the origin.  It is elementary to
check that $u$ constructed above is a solution
to $P_k u =0$ and the convergence of the integral defining $u$ and other standard
estimates follow in a manner analogous to that in \cite{HH1}, Lemma~2.1.  Using the
fact that the eigenfunction $f$ we have constructed is real-analytic at the
origin, we can easily show as in \cite{HH1}, Lemma~2.1 that the function $u$ is in
the Gevrey class
$2k+1$ at the origin. This Gevrey order agrees with the formula in \cite{Cha}  that
connects the location of the Tr\`eves strata in our example and the number of
commutation brackets one needs to descend to the center.\\

All these examples are   new. Of course, one can replace $ ( \left(\sum_{j=1}^p x_j^2\right)^{k}$ by a positive  elliptic 
 polynomial
 of order $m$ (with $m\geq 2 p$) in the variables $(x_1,\cdots,x_p)$. 

\appendix
\section{Schatten classes }\label{appendixA}
Here we collect a few well known results concerning Schatten classes.
 We refer to \cite{Sim} or \cite{BiSo} for more details.
We recall that a compact operator $A$ on an Hilbert space $\mathcal H$
 is in the Schatten class $\mathcal C _p$
 for some $p\in [1,+\infty[$ if the sequence $\mu_j$
 of the eigenvalues of $|A| =\sqrt{A^*A}$ satisfy
$\sum_j \mu_j^p < + \infty$.\\
When $p =1$, we speak about Trace class operators and, when $p=2$,  we recover
 the standard notion 
 of Hilbert-Schmidt operators.\\
When $p=1$, the trace map is defined by
\begin{equation}\label{deftr}
\mathcal C_1\ni A\mapsto \Tr A = \sum_j \langle A e_j\;|\; e_j \rangle\;,
\end{equation}
where $(e_j)$ is some orthonormal basis. It can be shown that this definition is independent of the choice
 of the basis and that the Trace map is continuous~:
\begin{equation}\label{conti}
| \Tr A \;| \leq   || \; |A| \; ||_{\mathcal C_1}\;.
\end{equation}
We have the H\"older relation, that is the
\begin{proposition}\label{PropCyc}.\\
 If $A \in \mathcal C_p$
 and $B \in \mathcal C_q$, then $AB \in \mathcal C _r$
 with $\frac 1r = \frac 1p + \frac 1q$.\\
Moreover, if $A\in \mathcal L(\mathcal H)$ and $B\in \mathcal C_q$,
 then $AB \in \mathcal C_q$.
\end{proposition}
When $r=1$, we will use constantly  the so-called cyclicity property of the trace~:
\begin{equation}\label{cyc}
\Tr (AB) = \Tr (BA)\;,\; \forall A \in \mathcal C_p, \forall B\in \mathcal C_q,
 \mbox{ with } \frac 1p + \frac 1q =1\;.
\end{equation}
The case $p=1$ is also true, if we replace $\mathcal C_\infty$ by $\mathcal L (\mathcal H)$. 
Various generalizations can be found in the book by M.~Birman and M.~Solomyak \cite{BiSo}.\\
Note also the property
\begin{equation}
|| A ||_{\mathcal C_1} = || A^* ||_{\mathcal C_1}\;.
\end{equation}
The following lemma will be useful for justifying extensions of the cyclicity
 rule.
\begin{lemma}.\\
We assume that $\mathcal H = L^2(\mathbb R^n)$.
Let $A$ be of class trace and $\chi$ a function in $C_0^\infty(\mathbb R^n)$ 
with compact support
 in a ball of radius $2$ and equal to $1$ on the ball of radius $1$.
Then if $A_j = \chi(\frac x j)A$ for $j \in \mathbb N^*$,  we have
 \begin{equation}\label{conv}
|| A - A_j  ||_{\mathcal C_1} \rightarrow 0\;,\;\mbox{  as }j\rightarrow + \infty\;;
\end{equation}
and 
\begin{equation}
\Tr A =  \lim_{j \rightarrow + \infty} A_j\;.
\end{equation}
\end{lemma}
{\bf Proof}.\\
Writing $A = |A|^\frac 12 C$ with $C$ Hilbert-Schmidt, we immediatly
 see that it is enough to treat the Hilbert-Schmidt case.
If one recalls that the Hilbert-Schmidt operators can be isometrically identified with the operators with
distribution  kernel in $L^2(\mathbb R^k \times \mathbb R^k)$,
 we are reduced to the application of the dominated convergence Theorem.
If $K$ is the kernel of $|A|^\frac 12$, we observe simply that
$$
\lim_{j \rightarrow + \infty} 
\int_{\mathbb R^n \times \mathbb R^n}
 (\chi (\frac{x}{j}) -1)^2 |K(x,y)|^2 dx dy = 0\;.
$$
We then conclude by observing that
$$
|| A - A_j||_{\mathcal C_1} \leq || (1 - \chi (\frac{\cdot}{j})) |A|^\frac 12 ||_{\mathcal C_2} \;\cdot \; || C ||_{\mathcal C_2}\;.
$$

\paragraph{Application}.\\
We use this lemma in the following context. We would like to show
 that
\begin{equation}\label{gencyc}
 \Tr (PC) = \Tr (CP) \;,
\end{equation}
where $P$ is a polynomial, $C$ is a trace class operator, such that
 $PC$ and $CP$ are trace class.
We first observe that the usual cyclicity trace rule gives~:
$$ \Tr (\chi (\frac{\cdot }{j})PC) = \Tr (CP\chi (\frac{\cdot }{j})) \;.$$
The lemma permits to justify the limiting procedure $j\rightarrow + \infty$.

Another trick could be to introduce an invertible operator $L$ such that
 $PL^{-1}$ is bounded and such that $L C$ is trace class.
Then one write~:
$$ \Tr  (PC) = \Tr (P L^{-1} L C) = \Tr ( L C P L^{-1}) $$
If $LCP$ and $ L^{-1}$ are in dual Schatten classes, one can reapply
 the cyclicity rule, and get
$$
 \Tr ( L C P L^{-1}) = \Tr ( L^{-1} L C P) = \Tr (CP)\;.
$$
All these conditions are practically easy to verify in the frame work
 of the pseudo-differential theory.

\section{Pseudodifferential operators and Schatten classes}\label{appendixB}
The theory of pseudo-differential operators gives an easy way for recognizing
that an operator belongs to a Schatten class. Let us recall a few elements of the theory. 
When $a$ belongs to a suitable class of symbols (see below),
 the Weyl quantization of the symbol $a$ consists in the introduction
 of the  operator 
$\mathcal S (\mathbb R^n)\ni u \mapsto  \Opw (a) u \in \mathcal S (\mathbb R^n)$ defined by~:
\begin{equation}\label{B.1}
 ( \Opw (a) u) (x) = (2 \pi)^{-n} \int \int \exp i <x-y, \xi> 
a(\frac{x+y}{2},\xi) u(y) dy d\xi\;.
\end{equation}
 As an extension
 of the Calderon-Vaillancourt theorem giving sufficient conditions
 for $L^2$-continuity, we have the following proposition for the Weyl-quantized pseudo-differential operators (See for example \cite{Ron}).
\begin{theorem}.\\
There exists $k$ depending only on the dimension such that, if
$$
N_{k,p} (a):= \sum_{|\alpha|\leq k } || D_{x,\xi}^\alpha a (x,\xi) ||_{L^p(\mathbb R^n \times \mathbb R^n)}  < + \infty
$$
then $\Opw (a)$ belongs to $\mathcal C_p$.
Moreover, we have for a suitable constant $C$:
\begin{equation}
|| \Opw (a) ||_{\mathcal C_p}\leq C\; N_{k,p} (a)\;.
\end{equation}
\end{theorem}
The Hilbert-Schmidt case (corresponding to $\mathcal C_2$)
 is more standard and we recall that~:
\begin{equation} \label{HS}
|| \Opw (a) ||_{\mathcal C_2}^2 = \int\int |a(x,\xi)|^2\; dx\, d\xi\;.
\end{equation}

The case $p= + \infty$ corresponds, when replacing $\mathcal C_\infty$ by $\mathcal L(L^2(\mathbb R^n))$, to the well known Calderon-Vaillancourt Theorem.

\section{On globally elliptic operators}\label{appendixC}
The last thing we would like to recall is the class of pseudodifferential operators adapted to our problem of analyzing the inverse of the operators
 $(- \Delta + P(x)^2)^{s}$. The reference \cite{HeRo}
 presents a pseudo-differential calculus which is exactly adapted
 to the situation. The symbols are indeed $C^\infty$ functions on 
$\mathbb R^n \times \mathbb R^n)$ for which there exists a real $M$ such that at $\infty$
\begin{equation}\label{C.1}
a(x,\xi) \sim \sum_{j\in \mathbb N}  a_{M-j} (x,\xi)\;,
\end{equation}
$a_{M-j}$ having the following homogeneity property
 for suitable $k >0$ and $\ell >0$
\begin{equation}\label{C.2}
a_{M-j} (\rho^k x, \rho^\ell \xi) = \rho^{M-j} a_{M-j} (x,\xi)\; ,\;\forall (x,\xi) \in \mathbb R^n\times \mathbb R^n\;, \forall \rho >0.
\end{equation}
We call this class $S^{M}_{k,\ell}$. We denote by $\Op S^{M}_{k,\ell}$
the class of operators defined as $\Opw  (a)$
 for some $a$ in $S^{M}_{k,\ell}$. We note that the composition
 of two operators $A_1 \in \Op S^{M_1}_{k,\ell}$ and
 of an operator  $A_2 \in \Op S^{M_1}_{k,\ell}$ gives 
$A_1 \circ A_2 \in  \Op S^{M_1+ M_2}_{k,\ell}$, the principal symbol of the product being
 simply the product
 of the principal symbols of $A_1$ and $A_2$.

The basic example is $L = - \Delta + P^2$ with $P$ homogeneous of degree $m$.
With $k = \frac 1m$, $\ell = 1$, we see that the symbol of this operator
 belongs to $S^{2}_{\frac 1m,1}$, so $L  \in \Op S^{2}_{\frac 1m,1}$ . This operator is
 ``elliptic'' in the sense that its principal symbol
 does not vanish on the sphere $S^{2n-1}$ and it is  shown in \cite{HeRo} that its
 inverse has a symbol in $S^{-2}_{\frac 1m,1}$.
Note also that a polynomial of order $k$ belongs to $S^{\frac {k}{m}}_{\frac 1m, 1} $. The question of determining if a pseudo-differential operator
 belongs to a Schatten class is then easy.
 The condition is simply
\begin{equation}\label{C.3}
\Opw (a) \in \mathcal C_p \mbox{ if } a \in S^{M}_{k,\ell} \mbox{ with }  M p + (k+ \ell) n < 0\;.
\end{equation}
\begin{remark}.\\
We note also that the pseudo-differential calculus gives an easy way for showing that the eigenvector whose existence is proved via Lidskii's  Theorem is actually in the Schwartz class $\mathcal S(\mathbb R^n)$.
\end{remark}

\end{document}